\documentclass[english,10pt]{article}
\usepackage[latin1]{inputenc}
\usepackage{amsmath,amsthm,amssymb}
\usepackage{graphicx}
\usepackage{color}

\newcommand{\be}{\begin{eqnarray}}
\newcommand{\ee}{\end{eqnarray}}

\def\eq#1{(\ref{#1})}

\def\neweq#1{\begin{equation}\label{#1}}
\def\endeq{\end{equation}}

\def\RR{{\mathbb R} }

\def\di{\displaystyle}
\def\ri{\rightarrow}
\def\intom{\int_\Omega}

\newtheorem{theorem}{Theorem}[section]

\title{\sc Existence and non--existence results\\ for quasilinear elliptic exterior problems\\ with nonlinear
boundary conditions
 }
\author{\sc  Roberta Filippucci, Patrizia Pucci and Vicen\c{t}iu R\u{a}dulescu}

\date{}

\begin{document}

\maketitle

\begin{abstract} Existence and non--existence results are established
for  quasilinear elliptic  problems with nonlinear
boundary conditions and lack of compactness. The proofs combine variational methods with
the geometrical feature,
due to the competition between the different growths of the non--linearities.
\end{abstract}

\section{Introduction and the main results}
Let $\Omega$ be a smooth exterior domain in $\RR^N$, that is,
$\Omega$ is the complement of a bounded domain with $C^{1,\delta}$
boundary ($0<\delta <1$). Assume that $p$ is a real number
satisfying $1<p<N$, $a\in L^\infty (\Omega)\cap
C^{0,\delta}(\overline\Omega)$ is a positive function, and $b\in
L^\infty (\Omega)\cap C(\Omega)$ is non--negative. Let
$p^*:=Np/(N-p)$ denote the critical Sobolev exponent. In Yu
\cite{yu} it is studied the following quasilinear problem
\neweq{yu1}\begin{cases}
\di -\mbox{div}\,
(a(x)|Du|^{p-2}Du)+b(x)|u|^{p-2}u=g(x)|u|^{r-2}u
\qquad\mbox{\rm in }\ \Omega,\\
\di u=0\quad \mbox{\rm on}\ \partial\Omega,\qquad
\di\lim_{|x|\ri\infty}u(x)=0,
\end{cases}
\endeq
where $p<r <p^*$ and $g\in L^\infty (\Omega)\cap L^{p_0}(\Omega)$,
with $p_0:=p^*/(p^*-r)$, is a non--trivial potential
which is positive on some non--empty open subset of $\Omega$.
Under these assumptions, Yu proved in \cite{yu} that  problem
\eq{yu1} has a weak positive solution $u$ of class
$C^{1,\alpha}(\overline\Omega\cap B_R(0))$ for any $R>0$ and some
$\alpha =\alpha(R)\in (0,1)$. Problems of this type are motivated
by mathematical physics (see, e.g., Reed and Simon \cite{reed} and
Strauss \cite{str}), where certain stationary waves in nonlinear
Klein--Gordon or Schr\"odinger equations can be reduced to this
form.

Actually, a weak solution  of \eqref{yu1} satisfies for all $\varphi\in E$ the identity
\be\label{sol}
\int_\Omega(a(x)|Du|^{p-2}Du\boldsymbol{\cdot}D\varphi+b(x)|u|^{p-2}u\varphi)dx=
\int_\Omega g(x)|u|^{r-2}u\varphi dx,\ee
where $E$ is the completion of $C^\infty_0(\Omega)$ under the underlying norm
$$\|u\|_{a,b}=\left(\int_\Omega[a(x)|Du|^p+b(x)|u|^p]dx\right)^{1/p}.$$
By Lemma~2 of \cite{yu} every weak solution $u$ of \eqref{yu1} is
in $L^q(\Omega)$ for every $q\in[p^*,\infty)$ and approaches 0 as $|x|\to\infty$.
Of course $E\sim H^{1,p}_0(\Omega)$ whenever $0<b_0\le b(x)\in L^\infty(\Omega)$.
Taking $\varphi=u$ in \eqref{sol} we get
$\|u\|_{a,b}^p=\|u\|^r_{L^r(\Omega;g)}$, so that \eqref{yu1} does not admit
nontrivial weak solutions whenever $g\le0$ a.e. in $\Omega$.

We consider a related problem involving
a mixed nonlinear boundary condition and  show that the above
result {\it does not} remain true in certain circumstances.
The main features of the paper are the following: (i) the quasilinear
differential operator $-\mbox{div}\,
(a(x)|Du|^{p-2}Du)$ in the left hand--side of \eqref{yu1} is affected by
a different perturbation which behaves like $|u|^{q-2}u$, where first $p<r<q<p^*$
and then $p<q<r<p^*$;
(ii) the Dirichlet boundary condition of \eqref{yu1} is replaced by a
mixed nonlinear boundary condition.
With the same hypotheses on $\Omega$, $a$, $g$, $p$ and $r$, we consider the problem
\neweq{1}
\left\{\begin{tabular}{lll} &$\di -\mbox{div}\,
(a(x)|Du|^{p-2}Du)+|u|^{q-2}u=\lambda g(x)|u|^{r-2}u$
\qquad & $\mbox{\rm in}\ \Omega,$\\
&$\di a(x)|Du|^{p-2}{\partial_\nu u}+b(x)|u|^{p-2}u=0$
\qquad &$\mbox{\rm on}\ \partial\Omega$,\\
\end{tabular} \right.
\endeq
where $\lambda$ is a real parameter and $\nu$ is the unit vector
of the outward normal on $\partial\Omega$. More precisely, we first assume
\begin{enumerate}
\item[(H1)] {\it $g\in L^\infty (\Omega)\cap L^{p_0}(\Omega)$,
with $p_0:=p^*/(p^*-r)$, $p<r<q<p^*$, is a non--negative function which is positive
on a non--empty open subset of $\Omega$};
\item[(H2)]
{\it $b$ is a continuous positive function on $\Gamma=\partial\Omega$.}
\end{enumerate}
Without altering the proof arguments below, the coefficient $1$ of the dominating term $|u|^{p-2}u$ can be replaced by any
function $f\in L^\infty(\Omega)$, with $\mbox{inf\,ess}_{\Omega}\,f>0$. Hence
equation \eqref{1} is the renormalized form.

Problem \eq{1} may be viewed as a prototype of pattern formation
in biology and is related to the steady--state problem for a
chemotactic aggregation model introduced by Keller and Segel
\cite{kels}. Problem \eq{1} also plays an important role in the
study of activator--inhibitor systems modeling biological pattern
formation, as proposed by Gierer and Meinhardt \cite{giem}.

By a {\it weak $($non--trivial$)$ solution} of problem \eqref{1} we mean a non--trivial
function $u\in X=E\cap L^q(\Omega)$ verifying
for all $\varphi\in X$ the identity
\be\label{sol1}
\int_\Omega a(x)|Du|^{p-2}Du\boldsymbol\cdot D\varphi dx+\int_\Gamma b(x)|u|^{p-2}u\varphi d\sigma&
+\displaystyle{\int_\Omega}|u|^{q-2}u\varphi dx
=\lambda\int_\Omega g(x)|u|^{r-2}u\varphi dx,\ee
where now $E$ is the completion of the restriction on $\Omega$ of functions of $C^\infty_0(\RR^N)$
with respect to the norm
$$\|u\|_{a,b}=\left(\int_\Omega a(x)|Du|^{p}dx+\int_\Gamma b(x)|u|^{p}d\sigma\right)^{1/p},$$
and $X$ is the reflexive Banach space endowed with the norm
$$\|u\|=\left\{\|u\|_{a,b}^p+\|u\|_{L^q(\Omega)}^{p}\right\}^{1/p}.$$
Hence, by (H1)--(H2), all the integrals in \eqref{sol1} are well
defined and converge.

The loss of compactness of the Sobolev imbeddings on unbounded
domains renders variational techniques more delicate. Some of the
papers treating problems on unbounded domains use special function
spaces where the compactness is preserved, such as spaces of radially symmetric
functions. We point out that even if $\Omega$ is unbounded, standard
compact imbeddings still remain true, e.g., if $\Omega$ is {\it
thin at infinity}, in the sense that $$\lim_{R\ri
\infty}\sup\left\{\mu (\Omega\cap B(x,1))\, :\ x\in\RR^N,\
|x|=R\right\}=0\,,$$ where $\mu$ denotes the Lebesgue measure and
$B(x,1)$ is the unit ball centered at $x$. Such arguments cannot
be applied to our general unbounded domain $\Omega$. In this case,
since $\Omega$ is not ``{\it thin}" and it looks like $\RR^N$ at
infinity (because $\Omega$ is an exterior domain), the analysis of
the compactness failure shows that a Palais--Smale sequence of the
associated energy functional (see Bahri and Lions \cite{bahlio})
differs from its weak limit by ``{\it waves}" that go to infinity.
However, the definition of $X$, combined with the main assumption
$p<r<p^*$, ensures that \be\label{cc}\mbox{\it the function
space $X$ is compactly embedded into the weighted Lebesgue space
$L^r\left(\Omega; g\right)$.}\ee

Taking $\varphi=u$ in \eq{sol1}, we have that any weak solution $u$ of \eq{1}
satisfies the equality
\neweq{2}\|u\|_{a,b}^p+\|u\|_{L^q(\Omega)}^q=\lambda\,\|u\|_{L^r(\Omega;g)}^r,
\end{equation}
so that
problem \eq{1} does not have any nontrivial solution whenever $\lambda\le0$.
 We first prove
that the result still remains true for  sufficiently small values of $\lambda>0$ when $p<r<q<p^*$,
that is,
the term $|u|^{q-2}u$ ``{\it dominates}" the right hand--side
and makes impossible the existence of a solution to our problem
\eq{1}. On the other hand, if $\lambda>0$ is sufficiently large,
then \eq{1} admits weak solutions. The precise statement of this
result is the following.

\begin{theorem}\label{t1} {\bf(The case $\boldsymbol{p<r<q<p^*}$).} Under the assumptions {\rm(H1)} and  {\rm(H2)}
there exists $\lambda^*>0$ such that
\begin{enumerate}
\item[\phantom{i}(i)] if $\lambda<\lambda^*$, then problem
\eqref{1} does not have any weak solution;
\item[(ii)] if $\lambda\ge \lambda^*$, then problem \eqref{1} has at least one weak solution $u$,
with the properties
\begin{enumerate}
\item[(a)] $u\in L^\infty_{\text{\scriptsize{\rm loc}}}(\Omega)$;
\item[(b)] $u\in C^{1,\alpha}(\Omega\cap B_R)$,
$\alpha=\alpha(R)\in(0,1)$;
\item[(c)] $u>0$ in $\Omega$;
\item[(d)] $u\in L^m(\Omega)$ for all $p^*\le m<\infty$ and $\lim_{|x|\ri\infty}u(x)=0$.
\end{enumerate}
\end{enumerate}
\end{theorem}

In the second part of the paper we consider condition (H1)$'$,
which is exactly assumption (H1), with the only exception that
condition $p<r<q<p^*$ is replaced by
$$p<q<r<p^*.$$

\begin{theorem}\label{t2}  Under the assumptions {\rm(H1)$'$} and  {\rm(H2)}
\begin{enumerate}
\item[\phantom{i}(i)]  problem \eqref{1} does not have any weak
solution for any $\lambda\le0$; \item[(ii)] problem \eqref{1} has
at least one weak solution $u$, with the properties (a)--(d) of
Theo\-rem~$\ref{t1}$ for all $\lambda>0$.
\end{enumerate}
\end{theorem}

\section{Proof of Theorem~\ref{t1}}
We point out in what follows the main ideas of the proof:

(a) There is some $\lambda^*>0$ such that problem \eqref{1} does
not have any solution for any $\lambda<\lambda^*$. This
means that if a solution exists then $\lambda$ must be
sufficiently large. One of the key arguments in this proof is based
on the assumption $q>r$. In particular, this proof yields an
energy lower bound of solutions in term of $\lambda$ which will be
useful to conclude that problem \eqref{1} has a non--trivial
solution if $\lambda=\lambda^*$.

(b) There exists $\lambda^{**}>0$ such that problem \eqref{1} has at least one
solution for any $\lambda>\lambda^{**}$. Next, by the properties of $\lambda^{*}$
and $\lambda^{**}$ we deduce that $\lambda^{**}=\lambda^{*}$. The proof uses variational arguments
and is based on the coercivity of the corresponding energy functional defined on $X$ by
$$J_\lambda(u)=\frac1p\| u\|^{p}_{a,b}
+\frac1q\|u\|_{L^q(\Omega)}^q
-\frac\lambda r\|u\|_{L^r(\Omega;g)}^r.$$
We show that the minimum of $J$ is achieved by a weak solution of \eq{1}.
In order to obtain that this global minimizer is not trivial, we prove
that the corresponding energy level is negative provided $\lambda$ is
sufficiently large.

{\it Step 1. Non--existence for $\lambda>0$ small.}
It is enough to show that, if there is a  weak solution of problem \eq{1},
then $\lambda$ must be sufficiently large. Assume that $u$ is a  weak
solution of \eq{1}, then  by \eqref{sol1} we get \eqref{2}.
Since $r<q$ and $g^{{q}/(q-r)}$
is in $L^1(\Omega)$ by (H1), applying the Young inequality we deduce
that
\neweq{3}
\lambda\|u\|_{L^r(\Omega;g)}^r\leq\frac{(q-r)\lambda^{{q}/(q-r)}}{q}\intom
g(x)^{{q}/(q-r)}dx+\frac rq\|u\|_{L^q(\Omega)}^q.
\end{equation}
Next, by \eq{2}, \eq{3}  and the fact that $u$ is non--trivial,
\neweq{4}\begin{aligned}
0<\| u\|_{a,b}^p&\leq\frac{q-r}{q}\,\lambda^{{q}/(q-r)}\intom g(x)^{{q}/{(q-r)}}dx
+\frac{r-q}{q}\|u\|_{L^q(\Omega)}^q\\
&\leq \frac{q-r}{q}\,\lambda^{{q}/(q-r)}\intom g(x)^{{q}/(q-r)}dx:= \lambda^{{q}/(q-r)}A<\infty.
\end{aligned}\end{equation}
The continuity of the imbedding $X\hookrightarrow L^r(\Omega
;g)$ implies that there exists $C=C(\Omega,g,p,q,r)>0$ such that
\neweq{5}
C\|v\|_{L^r(\Omega;g)}^p\leq\|v\|_{a,b}^p\end{equation}
for any $v\in X$. Thus, by \eq{2} and
\eq{5}, we have
$C\|u\|_{L^r(\Omega;g)}^p\leq
\lambda \|u\|_{L^r(\Omega;g)}^r.$
Since $p<r<q$, $\lambda>0$ and $\|u\|_{L^r(\Omega;g)}>0$ by \eqref{2},
we deduce that
$$\lambda\ge C\|u\|_{L^r(\Omega;g)}^{p-r}
\ge CC^{(-1+r/p)}\| u\|^{p-r}_{a,b}\ge C^{r/p}\lambda^{q(p-r)/p(q-r)}A^{(p-r)/p}.$$
It follows that
$\lambda\ge (A^{p-r}C^{\,r})^{(q-r)/r(q-p)},$ which also implies that $\lambda^*\le(A^{p-r}C^{\,r})^{(q-r)/r(q-p)}$.
This con\-clu\-des the proof of $(i)$.

In particular, Step 1. shows that if for some $\lambda>0$ problem
\eqref{1} has a  weak solution $u$, then
\begin{equation}\label{esten}
\left({C^{r}}/{\lambda^p}\right)^{1/(r-p)}\le\|u\|_{a,b}^{p}\le
\lambda^{q/(q-r)}A,\end{equation}
where $C=C(\Omega,g,p,q,r)>0$ is the constant given in \eqref{5}.
\smallskip

\noindent{\it Step 2. Coercivity of $J$.} It follows by (H1). Indeed, for any $u\in X$ and
all $\lambda>0$
$$
J_\lambda(u)=\frac1p\|u\|_{a,b}^p+\frac1{2q}\|u\|_{L^q(\Omega)}^q+\frac1{2q}\|u\|_{L^q(\Omega)}^q
-\frac\lambda r\|u\|_{L^r(\Omega;g)}^r.$$
By H\" older inequality and (H1) we have
\be\label{6}J_\lambda(u)\ge\frac1p\|u\|_{a,b}^p+\frac1{2q}\|u\|_{L^q(\Omega)}^q+\frac1{2q}\|u\|_{L^q(\Omega)}^q-
\frac\lambda r\|g\|_{L^{q/(q-r)}(\Omega)}\|u\|_{L^q(\Omega)}^r.
\ee
Now, since for any positive numbers $\alpha$, $\beta$, $q$ and $r$, with $r<q$, the function
$\Phi:\RR^+_0\to\RR$ defined by $\Phi(t)=\alpha t^r-\beta t^q$, achieves
its positive global maximum
$$\Phi(t_0)=\frac{q-r}q\left(\frac rq\right)^{r/(q-r)}
\alpha^{q/(q-r)}\beta^{r/(r-q)}>0$$ at point $t_0=\left({\alpha
r}/{\beta q}\right)^{1/(q-r)}>0,$ we immediately have $\alpha
t^r-\beta t^q\le C(q,r)\alpha^{q/(q-r)}\beta^{r/(r-q)}$, where
$C(q,r)=({q-r})\left(r^r/q^q\right)^{1/(q-r)}.$ Returning to
\eqref{6} and using the above inequality, with
$t=\|u\|_{L^q(\Omega)}$,
$\alpha=\lambda\|g\|_{L^{q/(q-r)}(\Omega)}/r$ and $\beta=1/2q$, we
deduce that
$$J_\lambda(u)\ge\frac1p\|u\|_{a,b}^p+\frac1{2q}\|u\|_{L^q(\Omega)}^q-
C(\lambda,q,r,g),$$
where $C(\lambda,q,r,g)=2^{r/(q-r)}(q-r)\left(\lambda\|g\|_{L^{q/(q-r)}(\Omega)}\right)^{q/(q-r)}/qr$.
This implies the claim.
\smallskip

Let $n\mapsto u_n$ be a minimizing sequence of $J_\lambda$ in $X$, which is bounded in $X$
by Step 2. Without
loss of generality, we may assume that $(u_n)_n$ is non--negative,
converges weakly to some $u$ in $X$ and converges also pointwise.
\smallskip

\noindent{\it Step 3. The  non--negative weak limit $u\in X$ is a weak solution of
\eqref{1}.} To prove this, we shall show that
$$J_\lambda(u)\le\liminf_{n\to\infty}J_\lambda(u_n).$$
By the weak lower semicontinuity of the norm $\|\cdot\|$ we have
$$\frac1p\|u\|_{a,b}^p+\frac1q\|u\|_{L^q(\Omega)}^q\le\liminf_{n\to\infty}
\left(\frac1p\|u_n\|_{a,b}^p+\frac1q\|u_n\|_{L^q(\Omega)}^q\right).$$
Next, the boundedness of $(u_n)_n$ in $X$ implies with the same
argument that
$$\|u\|_{L^r(\Omega;g)}=\lim_{n\to\infty}\|u_n\|_{L^r(\Omega;g)}$$
by \eqref{cc}. Hence $u$ is a global minimizer of $J_\lambda$ in $X$.
\smallskip

\noindent{\it Step 4. The weak limit $u$ is a non--negative  weak solution of
\eqref{1} if $\lambda>0$ is sufficiently large.} Clearly $J_\lambda(0)=0$. Thus, by Step 3
it is enough to show that there exists $\Lambda>0$ such that
$$\inf_{u\in X} J_\lambda(u)<0\quad\mbox{for all }\lambda>\Lambda.$$
Consider the constrained minimization problem \be\label{8}
\Lambda:=\displaystyle\inf\left\{\frac 1p\|w\|_{a,b}^p+\frac
1q\|w\|_{L^q(\Omega)}^q\,:\, w\in X\mbox{ and
}\|w\|_{L^r(\Omega;g)}^r=r\right\}.\ee Let $n\mapsto v_n\in X$ be
a minimizing sequence of \eqref{8}, which is clearly
bounded in $X$, so that we can assume, without loss of generality,
that it converges weakly to some $v\in X$, with
$\|v\|_{L^r(\Omega;g)}^r=r$ and
$$\di\Lambda=\frac1p\,\|v\|_{a,b}^p+\frac1q\,\|v\|_{L^q(\Omega)}^q$$
by the weak lower semicontinuity of $\|\cdot\|$.
Thus, $J_\lambda(v)=\Lambda-\lambda<0$ for any $\lambda>\Lambda$.
\smallskip

Now put
$$\begin{aligned}\lambda^{*}:&=\sup\{\lambda>0\,:\,\mbox{problem \eqref{1} does not admit any  weak
solution}\},\\
\lambda^{**}:&=\inf\{\lambda>0\,:\,\mbox{problem \eqref{1} admits a  weak
solution}\}. \end{aligned}$$
Of course $\Lambda\ge\lambda^{**}\ge\lambda^{*}>0$. To complete the proof of Theorem~\ref{t1}
it is enough to argue the following essential facts: (a) problem \eqref{1} has
a  weak solution for any $\lambda>\lambda^{**}$; (b) $\lambda^{**}=\lambda^{*}$
and problem \eqref{1} admits a  weak
solution when $\lambda=\lambda^{*}$.
\smallskip

\noindent{\it Step 5. Problem \eqref{1} has a  weak solution for
any $\lambda>\lambda^{**}$ and $\lambda^{**}=\lambda^{*}$.} Fix
$\lambda>\lambda^{**}$. By the definition of $\lambda^{**}$, there
exists $\mu\in(\lambda^{**},\lambda)$ such that that $J_\mu$ has a
non--trivial critical point $u_\mu\in X$. Of course, $u_\mu$ is a
sub--solution of \eqref{1}. In order to find a super--solution of \eqref{1}
which dominates $u_\mu$, we consider the constrained minimization
problem
$$\inf\left\{\frac1p\,\|w\|_{a,b}^p+\frac1q\,\|w\|_{L^q(\Omega)}^q-
\frac\lambda r\,\|w\|_{L^r(\Omega;g)}^r\,:\, w\in X\mbox{ and
}w\ge u_\mu\right\}.$$
Arguments similar to those used in Step 4
show that the above minimization problem has a solution
$u_\lambda\ge u_\mu$ which is also a weak solution of problem
\eqref{1}, provided $\lambda>\lambda^{**}$.

 We already know
that  $\lambda^{**}\ge\lambda^{*}$. But, by the definition of
$\lambda^{**}$ and the above remark,  problem \eqref{1} has no
solutions for any $\lambda <\lambda^{**}$. Passing to the supremum, this forces
$\lambda^{**}=\lambda^{*}$ and completes the proof.
\smallskip

\noindent{\it Step 6. Problem \eqref{1} admits a non--negative weak
solution when $\lambda=\lambda^*$.} Let $n\mapsto\lambda_n$ be a
decreasing sequence converging to $\lambda^{*}$ and let $n\mapsto
u_n$ be a corresponding sequence of non--negative weak solutions of
\eqref{1}. As noted in Step 2, the sequence $(u_n)_n$ is bounded
in $X$, so that, without loss of generality, we may assume that it
converges weakly in $X$, strongly in $L^r(\Omega ;g)$, and
pointwise to some $u^*\in X$, with $u^*\geq 0$. By \eqref{sol1},
for all $\varphi\in X$,
$$\int_\Omega a(x)|Du_n|^{p-2}Du_n\boldsymbol\cdot D\varphi dx+
\int_\Gamma b(x)|u_n|^{p-2}u_n\varphi d\sigma
+\displaystyle{\int_\Omega}|u_n|^{q-2}u_n\varphi dx
=\lambda_n\int_\Omega g(x)|u_n|^{r-2}u_n\varphi dx,$$
and passing to the limit as $n\to\infty$ we deduce that
$u^*$ verifies \eqref{sol1} for $\lambda=\lambda^*$, as claimed.

It remains to argue that $u^*\not=0$. A key ingredient in this
argument is the lower bound energy given in \eqref{esten}. Hence,
since $u_n$ is a non--trivial weak solution of problem \eq{1}
corresponding to $\lambda_n$, we have
$\|u_n\|_{a,b}^{p}\geq
\left({C^{r}}/{\lambda^p}\right)^{1/(r-p)}$ by
\eqref{esten}, where $C>0$ is the constant given in \eqref{5} and not depending
on $\lambda_n$. Next, since $\lambda_n\searrow\lambda^*$ as
$n\ri\infty$ and $\lambda^*>0$, it is enough to show that
\begin{equation}\label{ufinal}
\| u_n-u^*\|_{a,b}\ri 0\quad\mbox{as $n\ri\infty$}.\end{equation}

Since $u_n$ and $u^*$ are weak solutions
of \eqref{1} corresponding to $\lambda_n$ and $\lambda^*$, we have by \eqref{sol1},
with $\varphi=u_n-u^*$,
\begin{equation}\label{unustar}
\begin{aligned}
&\int_\Omega a(x)\left( |Du_n|^{p-2}Du_n
-|Du^*|^{p-2}Du^*\right)\boldsymbol\cdot D(u_n-u^*) dx\\
&\qquad+\int_\Gamma
b(x)\left(|u_n|^{p-2}u_n-|u^*|^{p-2}u^*\right)(u_n-u^*) d\sigma+
\int_\Omega (|u_n|^{q-2}u_n-|u^*|^{q-2}u^*)(u_n-u^*) dx\\
&=\int_\Omega
g(x)\left(\lambda_n\,|u_n|^{r-2}u_n-\lambda^*\,|u^*|^{r-2}u^*\right)(u_n-u^*)
dx.\end{aligned} \end{equation}
Elementary monotonicity properties imply that
$$\int_\Omega (|u_n|^{q-2}u_n-|u^*|^{q-2}u^*)(u_n-u^*) dx\ge 0\quad\mbox{and}\quad
\langle I'(u_n^*)-I'(u^*),u_n-u^*\rangle\ge 0,$$
where
$$I(u):=\|u\|_{a,b}^p/p.$$
 Since $\lambda_n\searrow\lambda^*$ as $n\ri\infty$ and $X$
is compactly embedded in $L^r(\Omega;g)$, for all $p>1$ relation \eq{unustar}
implies
\neweq{weakco}
0\le \langle
I'(u_n^*)-I'(u^*),u_n-u^*\rangle\le
\int_\Omega
g(x)\left[\lambda_n\,u_n^{r-1}-\lambda^*\,(u^*)^{r-1}\right](u_n-u^*)
dx \ri 0\endeq
as $n\ri\infty$.

Now, we distinguish the cases $p\geq 2$ and $1<p<2$
and we use the following elementary inequalities (see \cite[formula (2.2)]{si}): for all $\xi$, $\zeta\in\RR^N$
\neweq{diaz1}
|\xi-\zeta|^{p}\leq\begin{cases} c(|\xi|^{p-2}\xi- |\zeta|^{p-2}\zeta
)(\xi-\zeta)\qquad&\mbox{for }\phantom{1<\,}p\geq 2;\\
c\langle|\xi|^{p-2}\xi-|\eta|^{p-2}\eta,\xi-\eta \rangle^{p/2}
\left(|\xi|^p+|\eta|^p\right)^{(2-p)/2}\qquad&\mbox{for }1<p<2,\end{cases}\end{equation}
where $c$ is a positive constant.
\smallskip

\noindent
{\sc Case 1:} $p\geq 2$.
By \eq{diaz1} and \eq{weakco}, we immediately conclude that
$$\|u_n-u^*\|_{a,b}^p\leq c\langle I'(u_n^*)-I'(u^*),u_n-u^*\rangle
=o(1)\qquad\mbox{as $n\ri\infty$}.$$
\vskip 0.1cm

\noindent
{\sc Case 2:} $1<p< 2$. Since by convexity for all $\gamma\ge1$
\neweq{inega}(v+w)^\gamma\leq
2^{\gamma -1}(v^\gamma +w^\gamma)\,\qquad\mbox{for all $v,w\in
\RR^+_0$,}\endeq
then, for $\gamma=2/p$, we have
$$\| u_n-u^*\|^2_{a,b}\leq 2^{(2-p)/p}
\left[\left(\intom a(x)|D( u_{n}- u^*)|^{p}
\,dx\right)^{2/p}+\left(\int_\Gamma b(x)|u_{n}-u^*|^{p}
\,d\sigma\right)^{2/p}\right].
$$
Thus, in order to conclude that \eqref{ufinal} holds, it is enough to show that
$$\intom a(x)|D (u_{n}-
u^*)|^{p} \,dx\ri 0\qquad \mbox{and} \qquad\int_\Gamma
b(x)|u_{n}-u^*|^{p} \,d\sigma\ri 0$$
as $n\ri\infty$.
Indeed, combining \eq{diaz1} and \eq{inega}, we have
$$ \begin{aligned}&\int_\Omega a(x)|D (u_{n}- u^*)|^{p}\,dx\\
&\quad\leq c \int_\Omega a(x)\left\{(|D u_{n}|^{p-2}D u_{n}- |D
u^*|^{p-2}D u^*)\cdot D( u_{n}- u^*) \right\}^{p/2} \left(|D
u_{n}|^p+|D u^*|^{p}\right)^{(2-p)/2}\,dx\\
&\quad\leq c  \left(\int_\Omega a(x)(|D u_{n}|^{p-2}D u_{n}- |D
u^*|^{p-2}D u^*)\cdot D( u_{n}-u^*)\,dx
\right)^{p/2}\left(\|u_n\|_{a,b}^{p}+\|u^*\|_{a,b}^{p}\right)^{(2-p)/2}\\
&\quad\leq c  \left(\int_\Omega a(x)(|D u_{n}|^{p-2}D u_{n}- |D
u^*|^{p-2}D u^*)\cdot D( u_{n}- u^*)\,dx
\right)^{p/2}\left(\|u_{n}\|_{a,b}^{(2-p)p/2}+
\|u^*\|_{a,b}^{(2-p)p/2}\right)\\
&\quad \leq C_1 \left(\int_\Omega
a(x)(|D u_{n}|^{p-2}D u_{n}- |D u^*|^{p-2}D u^*)\cdot D( u_{n}-
u^*)\,dx \right)^{p/2},\end{aligned}$$
where $C_1=2c \left(\lambda^{q/(q-r)}A\right)^{(2-p)/2}$ by \eqref{esten} and $C_1$ is independent of $n$
by \eqref{4}.
Similar arguments yield
$$\int_\Gamma b(x)(u_{n}-u^*)^{p}
\,d\sigma\leq C_2\left(\int_\Gamma b(x)\left[ u_{n}^{p-1}-
(u^*)^{p-1} \right]( u_{n}- u^*)\,dx \right)^{p/2},$$
with an appropriate positive constant $C_2$ independent of $n$.
Combining
the above two inequalities with \eq{weakco} we conclude that
$\|u_n-u^*\|_{a,b}=o(1)$ as $n\ri\infty$, that is \eqref{ufinal} holds and $u^*$ is a
non--trivial non--negative weak solution of problem \eq{1}
corresponding to $\lambda=\lambda^*$.
\smallskip

Theorem~2.2 in Pucci and Servadei \cite{patraf1}, based on the Moser
iteration, shows that $u$ satisfies  {\it(a)}, since $u\in
W^{1,p}_{\text{\scriptsize{\rm loc}}}(\Omega)$, being $u\in X$,
$\bold A(x, u , \xi)=-a(x)| \xi|^{p-2}\xi$ and $B(x,u,\xi)=\lambda
g(x)|u|^{r-2}u-|u|^{q-2}u$ clearly verifies inequality (2.18) of
\cite{patraf1} by (H1); for other applications see also
\cite{patraf}. Next, again by the main assumptions on the
coefficient $a=a(x)$, an application of \cite[Corollary on
p.~830]{db}  due to DiBenedetto shows that the weak solution $u$
verifies also property {\it(b)}. Finally, {\it(c)} follows
immediately by the strong maximum principle since $u$ is a $C^1$
non--negative weak solution of the differential inequality
$\mbox{div}\, (a(x)|Du|^{p-2}Du)-|u|^{q-2}u\le0$ in $\Omega$, with
$q>p$, see, for instance, Section~4.8 of Pucci and Serrin \cite{ps2}
and the comments thereby. Property {\it(d)} follows using similar
arguments as in the proof of Lemma~2 of \cite{yu}, which is based on
Theorem~1 of Serrin \cite{s}.\qed

\section{Proof of Theorem~\ref{t2}}

Taking $\varphi=u$ in \eqref{sol1}, we see that any weak solution $u$ of \eqref{1}
satisfies the equality \eqref{2}, and the conclusion $(i)$ of Theorem~\ref{t2}
follows at once.
\smallskip

We next show that the $C^1$ energy functional $J_\lambda:X\to\RR$ satisfies the
assumptions of the Mountain Pass theorem of Ambrosetti and Rabinowitz
\cite{ar}. Fix $w\in X\setminus\{0\}$. Since $p<q<r$ then
$$J_\lambda(tw)=\frac{t^p}p\|w\|_{a,b}^p+\frac{t^q}{q}\|w\|_{L^q(\Omega)}^q
-\lambda\frac{t^r}r\|w\|_{L^r(\Omega;g)}^r<0$$
provided $t$ is sufficiently large. Next, by \eqref{cc}, \eqref{5}
and the fact that $p<q<r$ we observe that
$$J_\lambda(u)\ge\frac{1}q\|u\|^p-\frac{\lambda}r\|u\|_{L^r(\Omega;g)}^r\ge
\frac{1}q\|u\|^p-\frac{\lambda}{rC^{r/p}}\|u\|^r\ge\alpha>0,$$
whenever $\|u\|=\varrho$ and $\varrho>0$ is sufficiently small.  Set
$$\Gamma=\{\gamma\in C([0,1];X)\,:\,\gamma(0)=0,\,\,\gamma(1)\not=0\mbox{ and }
J_\lambda(\gamma(1))\le0\},$$
and put
$$c=\inf_{\gamma\in\Gamma}\max_{t\in[0,1]}J_\lambda(\gamma(t)).$$
Applying the Mountain Pass theorem without the Palais--Smale
condition we find a sequence $n\mapsto u_n\in X$ such that
\neweq{mplem1}J_\lambda(u_n)\to c\quad\mbox{and}\quad
J_\lambda'(u_n)\to0\end{equation} as $n\to\infty$. Moreover, since
$J_\lambda (|u|)\leq J_\lambda (u)$ for all $u\in X$, we can
assume that $u_n\geq 0$ for any $n\geq 1$. In what follows we
prove that $(u_n)_n$ is bounded in $X$. Indeed, since
$J_\lambda'(u_n)\to0$ in $X'$, then
$$\|u_n\|_{a,b}^p+\|u_n\|_{L^q(\Omega)}^q=\lambda\|u_n\|_{L^r(\Omega;g)}^r+o(1)$$
as $n\to\infty$. Therefore,
$$c+o(1)=J_\lambda(u_n)\ge\frac{1}q\|u_n\|^p-\frac{\lambda}r\|u_n\|_{L^r(\Omega;g)}^r\ge
\left(\frac1q-\frac1r\right)\left(\|u_n\|^p-1\right)+o(1)$$
as $n\to\infty$. Thus, since $q<r$, we deduce that the Palais--Smale sequence
 $(u_n)_n$ is bounded in $X$. Hence, up to a subsequence, we can assume that
 $(u_n)_n$  converges weakly in $X$ and strongly in ${L^r(\Omega;g)}$
to some element, say $u^*\geq 0$. From now on, with the same
arguments as in the proof of Theorem~\ref{t1}, we deduce that
$u^*$ is a weak solution of the problem \eq{1} such that
properties {\it(a)}--{\it(d)} are fulfilled. Due to the mountain-pass
geometry of our problem \eq{1} generated by the assumption
$p<q<r<p^*$, we are able to give the following alternative proof
in order to show that $u^*$ is a weak solution of \eq{1}. Fix
$\varphi\in C^\infty_0(\RR^N)$. Since $J_\lambda'(u_n)\to0$ in
$X'$, we have
$$\int_{\Omega}a(x)|Du_n|^{p-2}Du_n\cdot D\varphi dx+
\int_\Gamma b(x)u_n^{p-1}\varphi d\sigma + \int_{\Omega}u_n^{q-1}\varphi dx-\lambda
\int_{\Omega}g(x)u_n^{r-1}\varphi dx=o(1)$$ as
$n\ri\infty$.   Letting $n\ri\infty$, we deduce that
$$\int_{\Omega}a(x)|Du^*|^{p-2}Du^*\cdot D\varphi dx+
\int_\Gamma b(x)(u^*)^{p-1}\varphi d\sigma + \int_{\Omega}(u^*)^{q-1}\varphi dx-\lambda
\int_{\Omega}g(x)(u^*)^{r-1}\varphi dx=0$$
and so by density  $u^*$ satisfies
relation \eq{sol1}  for any $\varphi\in X$. It remains to
show that $u^*\not=0$. Indeed, by \eq{mplem1} and $n$ is sufficiently large we obtain
$$\begin{aligned} 0 \di < \frac c2&\leq J_\lambda (u_n)-\frac
1p\langle J_\lambda '(u_n),u_n\rangle\di \\&=\left(\frac 1q-\frac
1p\right)\|u_n\|_{L^q(\Omega)}^q-\lambda \left(\frac 1r-\frac
1p\right)\|u_n\|_{L^r(\Omega;g)}^r \leq\frac{\lambda (r-p)}{r}\|
u_n\|^r_{L^r(\Omega;g)},\end{aligned}$$ since $p<q<r$ . This implies
that $\| u^*\|^r_{L^r(\Omega;g)}>0$ and in turn $u^*\not=0$, as
required.

Finally, $u^*$ verifies properties {\it(a)}--{\it(d)}, as shown in the proof of Theorem~\ref{t1}. \qed
\medskip

\noindent {\bf Acknowledgments.} This work has been completed while
V.~R\u{a}dulescu was visiting the {\it Universit\`a degli Studi} of
Perugia in November 2006 with a GNAMPA--INdAM visiting professor
position. He is also supported by Grants 2-CEx06-11-18/2006 and
CNCSIS A-589/2007.  The first two authors were supported by the
Italian MIUR project titled {\it ``Metodi Va\-ria\-zio\-na\-li ed
Equazioni Differenziali non Lineari''}.

\bigskip\small

\noindent{\sc Dipartimento di Matematica e Informatica, Universit\`a degli Studi
di Perugia,  06123 Perugia, Italy}\\
{\it E-mail address}: {\tt roberta@dipmat.unipg.it}, {\tt pucci@dipmat.unipg.it}

\smallskip
\noindent{\sc Department of Mathematics, University of Craiova,
200585 Craiova, Romania}

\noindent{\sc Institute of Mathematics ``Simion Stoilow" of the
Romanian Academy, P.O. Box 1-764, 014700 Bucharest, Romania}\\
{\it E-mail address}: {\tt vicentiu.radulescu@math.cnrs.fr}

\end{document}